\documentclass{article}

\usepackage{latexsym,amssymb}
\title{A group theoretic condition equivalent to a condition on principal blocks}
\author{Geoffrey R. Robinson}

\begin{document}

\maketitle

\begin{abstract}
In this note, we consider the question: if $G$ is a finite group, which complex irreducible characters 
$\chi$ of $G$ belong to the principal $p$-block of $G$ for every prime divisor $p$ of $G$?  This question has been studied previously in a more general context by various authors, (see, for example Bessenrodt, Malle and Olsson,[1], and the references there). In particular, it is proved there that very few finite simple groups fail to have the property in question, and the exceptions are completely characterized.

\medskip
Here, we give a necessary and sufficient purely group-theoretic condition to ensure that the trivial character is the only irreducible character of $G$ which lies in the principal $p$-block for every prime divisor $p$ of $|G|$. To be precise, we prove a slightly more general statement about a set of distinct prime divisors of $|G|$ which is not necessarily exhaustive:

\medskip
\noindent {\bf Theorem 1.1:} \emph{ Let $G$ be a finite group, and let  $p_{1},p_{2}, \ldots, p_{n}$ be prime divisors of $|G|$. Then the following are equivalent:}

\medskip
\noindent i) \emph{ The trivial character is the only irreducible character of $G$ which lies in the principal $p_{i}$-block for $ 1 \leq i  \leq n.$}

\medskip
\noindent ii) \emph{The number of expressions of $g \in G$ as a product $g = x_{1}x_{2} \ldots x_{n}$ where each $x_{i} \in G$ is $p_{i}$-regular, is independent of the choice of $g$.}

\medskip
We also prove another variant  of this result which deals with $p$-sections other than the $p$-section of the identity.

\end{abstract}

\section{Introduction and notation}

\medskip
Let $G$ be a finite group. For $X$ a subset of $G,$ we let $X^{+} = \sum_{x \in X} x$ in the group algebra 
$\mathbb{C}G.$ For each prime divisor $p$ of $|G|$, we let $G_{p^{\prime}}$ denote the set of $p$-regular elements of $G$. For each complex irreducible character $\chi$ of $G,$ we let $\omega_{\chi}$ denote the linear character of $Z(\mathbb{C}G)$ associated to $\chi,$ which is defined by 
$$\omega_{\chi}(z) = \frac{\chi(z)}{\chi(1)}$$ for all $z \in Z(\mathbb{C}G).$ 

\medskip
If $\pi$ is a set of primes, and $n$ is a positive integer, then the integer $n_{\pi}$ is the largest divisor of $n$ which is only divisible by primes in $\pi.$

\medskip
Here we are interested in those finite groups $G$ which have the property that the only complex irreducible character of $G$ which lies in the principal $p$-block of $G$ for every prime divisor $p$ of $|G|$ is the trivial character. There are many finite simple groups with this property, (in fact, simple and quasi-simple groups with this property have been completely characterised by Bessenrodt, Malle and Olsson in the paper [1]). We will give here a necessary and sufficient  purely group-theoretic condition for the finite group $G$ to have a slightly more general character-theoretic property.

\medskip
We first prove a result which we view as somewhat analogous to a theorem of R. Brauer and H. Wielandt on a characterization of perfect groups via character theory (whose proof may be found, for example,  in Feit's book [3]):

\medskip
\noindent {\bf Theorem 1.1:} \emph{Let $G$ be a finite group, and let $p_{1},p_{2},\ldots ,p_{n}$ be distinct prime divisors of $|G|.$ Then the following are equivalent:}

\medskip
\noindent i) \emph{ The trivial character is the only irreducible character of $G$ which lies in the principal $p_{i}$-block for each $ i \leq n.$}

\medskip
\noindent ii) \emph{The number of expressions of $g \in G$ as a product $g = x_{1}x_{2} \ldots x_{n}$ where each $x_{i} \in G$ is $p_{i}$-regular, is independent of the choice of $g$.}

\medskip
\noindent{\bf Proof:} We note that (block orthogonality relations and) results of  Brauer-Feit [2], imply  that for each $i$, the irreducible character $\chi$ lies in the principal $p_{i}$-block of $G$ if and only if the linear character $\omega_{\chi}$ does not annihilate  $G_{p_{i}^{\prime}}^{+}$. Hence if only the trivial character lies in each principal $p_{i}$-block, it follows that  the product of all $n$ of the $G_{p_{i}^{\prime}}^{+}$'s is a rational multiple of the primitive idempotent of $Z(\mathbb{C}G)$ associated to the trivial character, so every element of $G$ appears with the same coefficient in the product, and then it is clear that the product must be a positive integer multiple of $G^{+}$.

\medskip
On the other hand, if the given product is an integer multiple of $G^{+},$ then each non-trivial irreducible character $\chi$ of $G$ annihilates the product, so the linear character $\omega_{\chi}$ must annihilate at least one factor $G_{p_{i}^{\prime}}^{+}.$ For any such $i,$ the character $\chi$ lies outside the principal $p_{i}$-block of $G$.

\medskip
\noindent {\bf Remark 1.2:} Letting $\pi = \{p_{1},p_{2}, \ldots ,p_{n} \}$, we note that if only the trivial character of $G$ lies in the principal $p_{i}$-block of $G$ for each $p_{i} \in \pi,$  then the number of times each $g$ is expressible in the form given in the theorem is a positive integer multiple of 
$|G|^{n-2} |G|_{\pi^{\prime}}.$ 

\medskip
For it was known to G. Frobenius that $|G_{p^{\prime}}|$ is always divisible by $|G|_{p^{\prime}}$ for any prime $p$. Hence the number of terms in the expansion of the product of the $|G_{p_{i}^{\prime}}|$'s is a positive integer multiple of $|G|^{n-1}|G|_{\pi^{\prime}},$ so since each element of $G$ occurs with equal multiplicity in the product, the claim follows. 

\section{Other $p$-sections}

\medskip
There are many variants of Theorem 1.1 available. For example, if we choose a prime $p$ and a $p$-element $z$ of $G$, we may consider the $p$-section $S_{z}$ of $z$ in $G$ (that is, the set of elements $g \in G$ whose $p$-part is $G$-conjugate to $z$). In the case that $z \in Z(P)$ for some Sylow $p$-subgroup $P$ of $G$, it is true (for the same reasons as before) that the irreducible character $\chi$ of $G$ belongs to the principal $p$-block of $G$ if and only if the linear character $\omega_{\chi}$ does not annihilate the $p$-section sum $S_{z}^{+}$ (this need no longer be true when $z$ is not central in any Sylow $p$-subgroup of $G$). Hence we may state a more general form of Theorem 1.1:

\medskip
\noindent {\bf Theorem 2.1:} \emph{Let $G$ be a finite group, and let $p_{1},p_{2},\ldots ,p_{n}$ be distinct prime divisors of $|G|.$ Let $P_{i}$ be a Sylow $p_{i}$-subgroup for each $i$, and for each $i$, choose $z_{i} \in Z(P_{i})$. Then the following are equivalent: }

\medskip
\noindent i) \emph{ The trivial character is the only irreducible character of $G$ which lies in the principal $p_{i}$-block for each $ i \leq n.$}

\medskip
\noindent ii) \emph{The number of expressions of $g \in G$ as a product $g = x_{1}x_{2} \ldots x_{n},$ where each $x_{i}$ lies in the $p$-section of $z_{i},$ is independent of the choice of $g$.}

\medskip
\noindent {\bf Acknowledgement:} We thank Gabriel Navarro for some helpful discussions related to this note.

\medskip
\section{References}

\medskip
\noindent [1] Bessenrodt, C.,  Malle, G. and Olsson, J.B.,  \emph{Separating characters by blocks}, Journal of the London Mathematical Society, (2),73,(2006),493-505.

\medskip
\noindent [2] Brauer, R. and Feit, W., \emph{On the number of irreducible characters of finite groups in a given block}, Proc Nat. Acad. Sciences USA,45,(1959),361-365.

\medskip
\noindent [3] Feit, W.,\emph{Characters of Finite Groups}, W.A. Benjamin,Inc., New York, (1967).

\end{document}